\documentclass {amsart}
\usepackage{amsmath,amssymb,amsbsy,amsfonts,amsthm,latexsym,amsopn,amstext,
                            amsxtra,euscript,amscd}

\newtheorem{conj}{Conjecture}

\newtheorem{problem}{Problem}

\def\bP{\mathbb P}
\def\bC{\mathbb C}
\def\C{\mathcal C}
\def\Z{\mathbb Z}
\def\X{\mathcal X}
\def\p{\mathfrak p}
\def\H{\mathcal H}
\def\M{\mathcal M}
\def\L{\mathcal L}
\def\iso{{\, \cong\, }}
\def\lar{\longrightarrow}
\def\bQ{{\mathbb Q}}
\def\bZ{{\mathbb Z}}

\def\p{\mathfrak p}
\def\X{\mathcal X}
\def\Aut{\mbox{Aut}}

\def\({\left(}
\def\){\right)}

\title{Computational algebra and algebraic curves}

\author{Tanush Shaska}

\begin{document}

\maketitle
\begin{center}Department of Mathematics, \\
University of Idaho,  \\
Moscow, ID, 83843
\end{center}

\begin{abstract}
The development of computational techniques in the last decade has made possible to attack some classical problems
of algebraic geometry.  In this survey, we briefly describe some open problems related to algebraic curves which
can be approached from a computational viewpoint.
\end{abstract}

\section{Introduction}
Computational algebra is a very active and rapidly growing field, with many applications to other areas of
mathematics, as well as computer science and engineering. In this survey, however, we will focus on applications
that computational algebra has on classical problems of mathematics and more explicitly on algebraic curves.

We survey topics related to automorphisms of algebraic curves, field of moduli versus field of definition, Hurwitz
groups and Hurwitz curves, genus 2 curves with split Jacobians etc. The problems we suggest are a very narrow
trend in algebraic geometry. However, they provide examples of how new computational techniques can be used to
answer some old questions.

In the second section we describe genus 2 curves with split Jacobians. There are many papers written on these
topic going back to Legendre and Jacobi in the context of elliptic integrals. The problem we suggest is to compute
the moduli space of  covers of degree 5, 7 from a genus 2 curve to an elliptic curve. This problem is completely
computational and could lead to some better understanding of some conjectures on elliptic curves; see Frey
\cite{Fr}.

In section three, we discuss the automorphism groups of algebraic curves. There has been some important progress
on this topic lately, however much more can be done. Extending some of the results to positive characteristic
would be important. Further we suggest computing the equations of Hurwitz curves of genus 14 and 17.

In section 4 we study hyperelliptic curves. Finding invariants which classify the isomorphism classes of
hyperelliptic curves of genus $g \geq 3$ is still an open problem. However, it is an easier problem to deal with
hyperelliptic curves with extra automorphisms. The main result here is from  \cite{GS} where dihedral invariants
were introduced which identify the isomorphism classes of such curves. Using these dihedral invariants one can
determine the automorphism group of hyperelliptic curves; see \cite{Sh5}. However, implementing such algorithm is
still a challenge since the loci of curves with prescribed automorphism group are still to be computed in terms of
dihedral invariants. The second problem of  section 4 is to find what solvable groups can occur as monodromy
groups of full moduli dimension for coverings of the Riemann sphere with a genus two curve; see section 4.1, for
details.

In the last section we focus on the field of moduli of algebraic curves. This is a classical problem of algebraic
geometry that goes back to Weil and  Shimura among many others.  An answer to the conjecture of section 5 would be
important in algebraic geometry, but also from a computational viewpoint. Problems 7 - 10  suggest some variations
of the field of moduli problem.

\medskip

{\bf Acknowledgment:} Most of the topics discussed in this paper are joint work with my collaborators. I would
like to thank J. Gutierrez, B. Guralnick, K. Magaard, S. Shpectorov, H. V\"olklein, M. Fried, J. Schicho, I.
Shparlinski among many others for many helpful discussions. This paper originated from my talk in {\sc ACA 03},
held in Raleigh, North Carolina. I want to thank the organizing committee of {\sc ACA 03}, especially M.
Giesbrecht , H. Hong, E. Kaltofen, and A. Szanto. Finally, I would like to thank E. Volcheck for suggesting that I
summarize my talk at ACA 03 in this article for the Communications of Computer Algebra.

\section{Genus 2 curves with split Jacobian}

First, we focus on genus 2 curves whose Jacobians are isogenous to a product of elliptic curves. These curves have
been studied extensively in the 19th century in the context of elliptic integrals. Legendre gave the first example
of such a curve and    then Jacobi, Clebsch, Hermite, Goursat, Brioschi, and Bolza explored them further. In the
late 20th century Frey and Kani, Kuhn,   Gaudry and Schost, Shaska and Voelklein, and many others have studied
these curves further. They are of interest for the  arithmetic    of genus 2 curves as well as elliptic curves.
See \cite{FK} for some conjectures that relate this topic with the arithmetic of  elliptic curves.

Let $C$ be a curve of genus 2 and $\psi_1:C \lar E_1$ a map of degree $n$, from $C$ to an elliptic curve $E_1$,
both curves defined    over  $\bC$. In \cite{Sh1}, we show that this map induces a degree $n$ map $\phi_1:\bP^1
\lar \bP^1$. We determine all possible ramifications  for $\phi_1$. If $\psi_1:C \lar E_1$ is maximal (i.e., does
not factor  non-trivially) then there exists a maximal map  $\psi_2:C\lar E_2$, of degree $n$, to some elliptic
curve $E_2$ such that there is an isogeny of degree $n^2$ from the Jacobian $J_C$  to $E_1 \times E_2$. We say
that $J_C$ is  $(n,n)$-decomposable. If the degree $n$ is odd the pair $(\psi_2, E_2)$ is  canonically determined;
see \cite{Sh1} for details.

We denote the moduli space of such degree $n$ coverings $\phi: \bP^1\to \bP^1$ by $\L_n$. This space is studied by
Kani and it is called ``modular diagonal space''.  It can be viewed also as the Hurwitz space of covers $\phi:
\bP^1 \to \bP^1$ with ramification determined above. For our purposes, $\L_n$ will simply be the locus of genus 2
curves whose Jacobian is $(n,n)$-isogenous to a product of two elliptic curves.

The locus $\L_2$ of these genus 2 curves  is a 2-dimensional subvariety of the moduli space $\mathcal M_2$  and
is studied in detail in \cite{ShV1}. An equation for $\L_2$ is already in the work of Clebsch and Bolza. We use a
birational parametrization of $\L_2$ by affine 2-space to study the relation between the j-invariants of the
degree 2  elliptic subfields. This extends work of Geyer, Gaudry, Stichtenoth and others. We find a 1-dimensional
family of genus  2 curves having exactly two isomorphic elliptic subfields of degree 2; this family is
parameterized by the  j-invariant of these subfields. This was a joint project with H. V\"olklein, published in
the proceedings of  professor Abhyankar's 70th birthday conference.

If $n >2 $, the surface $\L_n$ is less understood. The case $n=3$ was initially studied by Kuhn \cite{Ku} where
some computations for $n=3$ were performed. The computation of the equation of $\L_3$ was a major computational
effort. A detailed description of this computation is given in \cite{Sh2}. Computational algebra techniques (i.e.,
Groebner basis, Buchberger algorithm) and computational algebra packages (i.e, Magma, Maple, GAP) were used. In
\cite{Sh2}, we study genus 2 function fields $K$ with degree 3 elliptic subfields. We show that the number of
$Aut(K)$-classes of such subfields of fixed $K$ is 0,1,2 or 4. Also we compute an equation for the locus of such
$K$ in the moduli space of genus 2 curves.  Equations of $\L_n$ for $n > 5$ are still unknown.

Let $\C$ be a genus 2 curve defined over $k$, $char (k) =0$. If ${\bar k} (\C) $ has a degree 3 elliptic subfield
then   the automorphism group $Aut(\C)$ is isomorphic to one of the following: $\bZ_2, V_4, D_8$, or $D_{12}$,
where  $D_n$ is the dihedral group of order $n$. There are exactly six genus two curves $\C$ defined over $\bC$
with $Aut(\C)$  isomorphic to $D_8$ (resp., $D_{12}$). We further show that only four (resp., three) of the curves
with group $D_8$ (resp., $D_{12}$) are defined over $\bQ$. This is summarized in \cite{Sh3}.

Continuing on the work of the above papers, we suggest the following problem:

\begin{problem}
Determine the locus $\L_n$ in $\M_2$ for $n= 5, 7$. Further, determine the relation between the elliptic curves
$E_1$ and $E_2$ in each case.
\end{problem}
Using techniques from \cite{ShV1, Sh2} this becomes simple a computational problem. However, determining such loci
requires the use of a Groebner basis algorithm. Computationally this seems to be difficult for $n=5, 7$.

\section{The automorphism group of a compact Riemann surface}
%
Computation of automorphism groups of compact Riemann surfaces is a classical problem that goes back to Schwartz,
Hurwitz, Klein, Wiman and many others. Hurwitz showed that the order of the automorphism group of a compact
Riemann surface of genus $g$ is at most $84 (g-1)$, which is known as the Hurwitz bound. Klein was mostly
interested with the real counterpart of the problem, hence the term   ``compact Klein surfaces''. Wiman studied
automorphism groups of hyperelliptic curves and orders of single automorphisms.

The 20th century produced a huge amount of literature on the subject.  Baily \cite{Ba} gave an analytical proof of
a theorem of Hurwitz: if  $g \geq 2$, there exists a curve of genus $g$ with non-trivial   automorphisms. In other
papers was treated the number of  automorphisms of a Riemann surface; see Accola \cite{Ac1}, Maclachlan
\cite{Mc1}, \cite{Mc2} among others. Accola \cite{Ac2} gives a  formula relating the genus of a Riemann surface
with the subgroups of  the automorphism group; known as Accola's theorem. Harvey studied  cyclic groups and Lehner
and Newman maximal groups that occur as  automorphism groups of Riemann surfaces.

A group of automorphisms of a compact Riemann surface $X$ of genus $g$  can be faithfully represented via its
action on the abelian   differentials on $X$ as a subgroup of $GL(g,\bC)$. There were many efforts to classify the
subgroups $G$ of $GL (g,\bC)$ that so arise, via the cyclic subgroups of $G$ and conditions on the matrix elements
of $G$. In a series of papers, I. Kuribayashi, A. Kuribayashi, and Kimura compute the lists of subgroups which
arise this way for $g=3,4$, and 5.

By covering space theory, a finite group $G$ acts (faithfully) on a genus $g$ curve if and only if it has a genus
$g$ generating system; see \cite{MS}. Using this purely group-theoretic condition, Breuer \cite{Br} classified all
groups that act on a curve of genus $\le 48$. This was a major computational effort using the computer algebra
system GAP. It greatly improved on several papers dealing with small genus, by various authors.

Of course, for each group in Breuer's list, all subgroups are also in the list. This raises the question how to
pick out those groups that occur as the {\bf full automorphism group} of a genus $g$ curve. This question is
answered in the following paper.

Let $G$ be a finite group, and $g\ge2$. In a joint project with Magaard, Shpectorov, and V\"olklein we study the
locus of genus $g$ curves that admit a $G$-action of given type, and inclusions between such loci; see \cite{MS}.
We use this to study the locus of genus $g$ curves with prescribed automorphism group $G$. We completely classify
these loci for $g=3$ (including equations for the corresponding curves), and for $g\le10$ we classify those loci
corresponding to ``large'' $G$.

We suggest the following:
\begin{problem}
Determine the list of possible automorphism groups of algebraic curves of small genus (i.e., $g \leq 10$) in every
characteristic.
\end{problem}
For $g=2$ this list is well known (it appears also in \cite{ShV1}). For $g=3$ it can probably be completed from
work of Brock, Wolper and others. However, for $g > 3$ such list of groups is unknown. It would be nice to have a
complete list for ``small genus'', say $g \leq 10$. Since, such lists tend to grow as genus grows,  such
information could be organized in a database and be very helpful to the mathematics community.  It is important to
mention that such lists are not known even in characteristic zero.

In \cite{Sh5} (ISSAC 03)  a new algorithm was introduced to compute the automorphism group of a given
hyperelliptic curve. However, this  will be discussed in more detail in the next section.

\subsection{Hurwitz curves}

A Hurwitz curve is a genus $g$ curve, defined over an algebraically closed field of characteristic zero, which has
$84(g-1)$ automorphisms. A group $G$ that can be realized as an automorphism group of a Hurwitz curve is called a
Hurwitz group. There are a lot of papers by group-theoretists on Hurwitz groups, surveyed by Conder. It follows
from Hurwitz's presentation that a Hurwitz group is perfect. Thus every quotient is again a Hurwitz group, and if
such a quotient is minimal then it is a non-abelian simple group. Several infinite series of simple Hurwitz groups
have been found by Conder, Malle, Kuribayashi, Zalessky, Zimmermann and others. In 2001, Wilson showed the monster
is a Hurwitz group; see \cite{MS} for a complete list of references.

Klein's quartic is the only Hurwitz curve of genus $g\le3$. Fricke showed that the next Hurwitz group occurs for
$g=7$ and has order 504. Its group is $SL(2,8)$, and an equation for it was computed by Macbeath in 1965. Klein's
quartic and Macbeath's curve are the only Hurwitz curves whose equations are known. Further Hurwitz curves occur
for $g=14$ and $g=17$ (and for no other values of $g\le 19$). It is natural, to try to write equations for these
Hurwitz curves of genus 14, 17.

\begin{problem}
Compute equations for the Hurwitz curves of genus 14, and possibly 17.
\end{problem}

\section{Computational aspects of hyperelliptic curves}

It is an interesting and difficult problem in algebraic geometry is to obtain a generalization of the theory of
elliptic modular functions to the case of higher genus. In the elliptic case this is done by the so-called
$j$-{\it invariant} of elliptic curves. In the case of genus $g=2$, Igusa (1960) gives a complete solution via
{\it absolute invariants} $i_1, i_2, i_3$ of genus 2 curves. Generalizing such results to higher genus is much
more difficult due to the existence of non-hyperelliptic curves. However, even restricted to the hyperelliptic
moduli $\H_g$ the problem is still unsolved for $g \geq 3$. In other words, there is no known way of identifying
isomorphism classes of hyperelliptic curves of genus $g\geq 3$. In terms of classical invariant theory this means
that the field of invariants of binary forms of degree $2g+2$ is not known for $g\geq 3$.

The following is a special case of $g=3$.
\begin{problem}
Find invariants which classify the isomorphism classes of genus 3 hyperelliptic curves.
\end{problem}
This is equivalent with determining the field of invariants of binary octavics. The covariants of binary octavics
were determined in 1880 by von Gall. The generators of the ring of invariants were determined by Shioda in 1965.
However, the field of invariants is unknown. This is a computational problem and should be possible to solve with
now-day techniques.  Extending this to positive characteristic would be quite interesting.

In a joint paper with J. Gutierrez, we find invariants that identify isomorphism classes of genus $g$
hyperelliptic curves with extra (non-hyperelliptic) involutions; see \cite{GS}. This result gives a nice way of
doing computations with these curves. We call such invariants {\it dihedral invariants} of hyperelliptic curves.
Let $\L_g$ be the locus in $\H_g$ of hyperelliptic curves with extra involutions. $\L_g$ is a $g$-dimensional
subvariety of $\H_g$. The dihedral invariants yield a birational parametrization of $\L_g$. Computationally these
invariants give an efficient way of determining a point of the moduli space $\L_g$. Moreover, we show that the
field of moduli is a field of definition (see below) for all $\p \in \L_3$ such that $|\Aut (\p)|>4$.

Dihedral invariants can be used to study the field of moduli of  hyperelliptic curves in $\L_g$ (cf. section 5).
Whether or not the field of moduli is a field of definition is in general a difficult problem that goes back to
Weil, Shimura et al. In ASCM 2003, we conjecture that for each $\p \in \H_g$ such that $|\Aut (\p)| > 2$ the field
of moduli is a field of definition. Making use of dihedral invariants we show that if the Klein 4-group can be
embedded in the reduced automorphism group of $\p\in \L_g$ the our conjecture holds; see \cite{Sh4} for details.

The families of hyperelliptic curves with reduced automorphism group (i.e., the automorphism group modulo the
hyperelliptic involution) isomorphic to $A_4$ or a cyclic group, are studied in \cite{Sh6}. We characterize such
curves in terms of classical invariants of binary forms and in terms of dihedral invariants. Further, we describe
algebraically the loci of such curves for $g\leq 12$ and show that for all curves in these loci the field of
moduli is a field of definition.

New techniques for computing the automorphism group of a genus $g$ hyperelliptic curve $\X_g$ are discussed in
\cite{Sh5}. The first technique uses the classical $GL_2 (k)$-invariants of binary forms. This is a practical
method for curves of small genus, but has limitations as the genus increases, due to the fact that such invariants
are not known for large genus. The second approach, which uses dihedral invariants of hyperelliptic curves, is a
very convenient method and works well in all genera. We define the normal decomposition of a hyperelliptic curve
with extra automorphisms. Then, dihedral invariants are defined in terms of the coefficients of this normal
decomposition. We define such invariants independently of the automorphism group $\Aut (\X_g)$. However, to
compute such invariants the curve is required to be in its normal form. This requires solving a nonlinear system
of equations. We discover conditions in terms of classical invariants of binary forms for a curve to have reduced
automorphism group $A_4$, $S_4$, $A_5$.

In the case of hyperelliptic curves the list of groups are completely determined in characteristic zero by work of
Bujalance, Gromadzky, and Gamboa. We suggest the following:
\begin{problem}
Implement a fast algorithm that does the following: Given a genus $g $ hyperelliptic curve $\X_g$, determine the
automorphism group of $\X_g$.
\end{problem}
The known algorithms (even the recent ones) approach the problem by solving a system of equations via Groebner
basis. This is normally inefficient and expensive.  We have implemented such programs for small $g$ and these
results can be extended even further. It will be valuable to organize such results in a computer algebra package
and extend to $g \leq 10$.

\subsection{The monodromy group of a genus 2 curve covering $\mathbb P^1$}
\def\sem{{\rtimes}}

Determining the monodromy group of a generic genus $g$ curve covering $\mathbb P^1$ is a problem with a long
history which goes back to Zariski and relates to Brill-Nother theory. Let $\X_g$ be generic curve of genus $g$
and $f: X_g \to {\mathbb P}^1$ a degree $n$ cover. Denote by $G:=Mon (f)$, the monodromy group of $f: X_g \to
{\mathbb P}^1$. Zariski showed that for $g > 6$, $G$ is not solvable. For $ g \leq 6$ the situation is more
technical. This has been studied by many authors e.g., Fried, Guralnick, Neubauer, Magaard, V\"olklein et al.
However, the problem is open for $g=2$. Guralnick and Fried, in a preprint dated at 1986, have shown that for $G$
primitive in $S_n$ and solvable there are six possibilities for $G$. Two of those are obvious cases $S_3, S_4$.
The other four groups are  $ D_{10}, \, \Z_3^2 \sem D_8, \, AGL_2 (3), \, S_4 \wr \Z_2 $; see \cite{FG}. The
corresponding signatures are:
$$(2^2, 2^2, 2^2, 2^2, 2^2, 2^2), \, \, (2^3, 2^3, 2^3, 2^3, 2^4, 2^4), $$
$$(2^3, 2^3, 2^3, 2^3, 3^2, 3^2), \, \, (2^6, 2^6, 2^6, 2^4, 2^4, 2^4, 2^4).
$$
\begin{problem}
For each case above, determine the locus of such genus 2 curves in $\M_2$ (e.g., the equation of such locus in
terms of invariants $i_1, i_2, i_3$) and its dimension.
\end{problem}
Notice that via the braid group action (using GAP), we can show that the  corresponding Hurwitz spaces are
irreducible. We expect in all cases that   the dimension of the locus in $\M_2$ is $\leq 2$.

\section{Field of moduli versus the field of definition}
Let $\X$ be a curve defined over $k$. A field $F \subset k$ is called a {\it field of definition} of $\X$ if there
exists $\X'$ defined over $F$ such that $\X \iso \X'$. The {\bf field of moduli} of $\X$ is a subfield $F \subset
k$ such that for every automorphism $\sigma$ of $k$ $\X$ is isomorphic to $\X^\sigma$ if and only if $ \sigma_F =
id$.

The field of moduli is not necessary a field of definition. To determine the points $\p \in \M_g$ where the field
of moduli is not a field of definition is a classical problem in algebraic geometry and has been the focus of many
authors, Weil, Shimura, Belyi, Coombes-Harbater, Fried, D\'ebes, Wolfart among others.

Weil (1954) showed that for every algebraic curve with trivial automorphism group, the field of moduli is a field
of definition. Shimura (1972) gave the first example of a family of curves such that the field of moduli is not a
field of definition. Shimura's family were a family of hyperelliptic curves.  Further he adds: {\it `` ... the
above results combined together seem to indicate a rather complicated nature of the problem, which almost defies
conjecture. A new viewpoint is certainly necessary to understand the whole situation''}

It seems as hyperelliptic curves provide the most interesting examples. For example, we are not aware of any
explicit examples of non-hyperelliptic curves such that the field of moduli is not a field of definition.
Moreover, with the help of dihedral invariants we have a way of describing the points of moduli in the locus
$\L_g$. Hence, we focus on hyperelliptic curves.

We call a point $\p \in \H_g$ a {\it moduli point}. The field of moduli of $\p$ is denoted by $F_\p$. If there is
a curve $\X_g$ defined over $F_\p$ such that $\p=[\X_g]$, then we call such a curve a {\it rational model over the
field of moduli}. Consider the following problem:

\smallskip

{\it Let the moduli point $\p \in \H_g$ be given. Find necessary and sufficient conditions that the field of
moduli $F_\p$ is a field of definition. If $\p$ has a rational model $\X_g$ over its field of moduli, then
determine explicitly the equation of $\X_g$.}

\smallskip

In 1993, Mestre solved the above problem for genus two curves with automorphism group $\Z_2$. Mestre's approach is
followed by Cardona and Quer (2002) to prove that for points $\p \in \M_2$ such that $| Aut (\p)| > 2$ the field
of moduli is a field of definition; see also \cite{Sh3} for a different approach. In his talk at ANTS V (see
\cite{Sh3}), the author conjectured the following:
\begin{conj} Let $\p \in \H_g $ be a moduli point such that $|Aut(\p)| >
2$. Then, its field of moduli  is a field of definition.
\end{conj}
The author has proved this conjecture for curves with reduced automorphism group isomorphic to $A_4$ and genus $g
\leq 12$; see \cite{Sh6}. Also the conjecture is true for $g=3$ and $| \Aut (\p) | > 4$; see \cite{GS}.
Furthermore, we intend to investigate the conjecture in all cases:
\begin{problem}
Investigate Conjecture 1 in all cases.
\end{problem}
In studying the above conjecture, we are looking for more than just a  true or false answer. We would like a way
to determine the field of moduli of any hyperelliptic curves with extra automorphisms. Generically, dihedral
invariants accomplish this for curves with extra involutions (i.e., locus $\L_g$). However, there is also the
singular locus in $\L_g$ which needs to be considered. And then, there are also hyperelliptic curves with extra
automorphisms which are not in $\L_g$. The upshot would be to solve the following:
\begin{problem}
Let $\p \in \H_g $. Determine if the field of moduli is a field of definition. In that case, explicitly find a
rational model of the curve over its field of moduli.
\end{problem}
The above problems lead to the following:
\begin{problem}
Find necessary and sufficient conditions in terms of invariants of binary forms such that a hyperelliptic curve
has no extra automorphisms.
\end{problem}
Such conditions were known to Clebsch and Bolza for $g=2$. These conditions were used by Mestre in \cite{Me}.
Finding similar conditions for $g\geq 3$ would help extend Mestre's algorithm to $g \geq 3$. Solving the above
problem would give a way of investigating Conjecture 1 without the hypothesis $| \Aut (\p) | >2$.
\begin{problem}
Find an algorithm which does the following: Let $\p \in \H_g $ such that $|\Aut (\p)| = 2 $. Determine if the
field of moduli is a field of definition.
\end{problem}



\begin{thebibliography} {99}

\bibitem[Ac1] {Ac1} {\sc R. Accola}, On the number of automorphisms of a closed Riemann surface,
Trans. Amer. Math. Soc. {\bf 131} (1968), 398-408.

\bibitem[Ac2] {Ac2} {\sc R. Accola},
Two theorems on Riemann surfaces with noncyclic automorphism groups, Proc. Amer. math. Soc. {\bf 25} (1970),
598-602.

\bibitem[Ba]{Ba} {\sc W. Baily},
On the automorphism group of a generic curve of genus $>2$, {\it J. Math. Kyoto Univ.} {\bf 1}
(1961/1962), 101--108; correction, 325.

\bibitem[Br] {Br} {\sc Th. Breuer},
Characters and automorphism groups of compact Riemann surfaces, London Math. Soc. Lect. Notes {\bf 280}, Cambridge
Univ. Press 2000.

\bibitem[Bo] {Bo}
{\sc O. Bolza}, On binary sextics with linear transformations into themselves, Amer. J. Math. 10 (1888), 47-70.

\bibitem[BS] {BS} {\sc R. Brandt and H. Stichtenoch},
Die Automorphismengrupenn hyperelliptischer Kurven. {\it Manuscripta Math} {\bf 55} (1986), no. 1, 83--92.


\bibitem[Bu] {Bu} {\sc E. Bujalance E, J. M. Gamboa, G. Gromadzki},
The full automorphism groups of hyperelliptic Riemann surfaces, {\it Manuscripta Math.} {\bf 79} (1993), no. 3-4,
267--282.


\bibitem[CF] {CF} {\sc J. W. S. Cassels and E. V. Flynn};
Prolegomena to a Middlebrow Arithmetic of Curves of Genus Two, LMS, 230, 1996.

\bibitem[Cl] {Cl}
{\sc A. Clebsch}, Theorie der Bin{$\ddot{a}$}ren Algebraischen Formen, Verlag von B.G. Teubner, Leipzig (1872).


\bibitem[DM] {DM} {\sc P. Deligne P, D. Mumford D},
The irreducibility of the space of curves of given genus, Publ. Math.Hautes \'Etudes Sci. {\bf 36}, 75-109, 1969.

\bibitem[ES] {ES}
{\sc T. Ekedahl T, J. P. Serre}, Exemples de courbes alg\'ebriques \'a jacobienne compl\'etement d\'ecomposable.
{\it C. R. Acad. Sci. Paris SÈr. I Math.}, 317 (1993), no. 5, 509--513.


\bibitem[Fr] {Fr} {\sc G. Frey}, On elliptic curves with isomorphic torsion
structures and corresponding curves of genus 2. {\it Elliptic curves, modular forms, and Fermat's last theorem
(Hong Kong, 1993)}, 79-98, Ser. Number Theory, I, {\it Internat. Press, Cambridge, MA}, (1995).

\bibitem[FK] {FK} {\sc G. Frey and E. Kani}, Curves of genus 2 covering elliptic
curves and an arithmetic application. {\it Arithmetic algebraic geometry (Texel, 1989)}, 153-176, {\it Progr.
Math.}, 89, {\it Birkh\"auser Boston, Boston, MA, (1991)}.

\bibitem [FG] {FG} {\sc M. Fried and R. Guralnick},
Radicals don't uniformize the generic curve of genus $g > 6$. (preprint).

\bibitem[FMV] {FMV} {\sc G. Frey, K. Magaard, and H. V\"olklein}
The general curve covers $\mathbb P$ with monodromy group $A_n$, preprint.

\bibitem[Gu] {Gu6} {\sc J. Gutierrez},
 A polynomial decomposition algorithm
over factorial domains, Comptes Rendues Mathematiques, de Ac. de Sciences, 13 (1991), 81-86.


\bibitem[Gu] {Gu2}
{\sc J. Gutierrez, H. Niederreiter, I. Shparlinski}, On the multidimensional distribution of inversive
congruential pseudorandom numbers in parts of the period. Monatsh. Math. 129 (2000), no. 1, 31--36.

\bibitem[GS] {GS}
{\sc J. Gutierrez and T. Shaska}, Hyperelliptic curves with extra involutions, (submitted), 2002.

\bibitem[Ku] {Ku} {\sc M. R. Kuhn}, Curves of genus 2 with split Jacobian. {\it Trans. Amer. Math. Soc}
{\bf 307}, 41-49, 1988.

\bibitem[Ig] {Ig} {\sc J. Igusa}, Arithmetic Variety Moduli for genus 2. {\it
Ann. of Math}. (2), 72, 612-649, 1960.


\bibitem[Mc1] {Mc1}
{\sc C. Maclachlan}, Abelian groups of
      automorphisms of compact Riemann surfaces, {\it Proc. London Math.
      Soc. (3)} {\bf 15} (1965), 699--712.

\bibitem[Mc2] {Mc2} {\sc C. Maclachlan}, A bound for the number of
      automorphisms of a compact Riemann surface, {\it J. London Math. Soc.
      (2)} {\bf 44} (1969), 265--272.


\bibitem[MS] {MS} {\sc K. Magaard, T. Shaska, S. Shpectorov, and H.
V\"olklein}, The locus of curves with prescribed automorphism group. Communications in arithmetic fundamental
groups (Kyoto, 1999/2001). S\=urikaisekikenky\=usho K\=oky\=uroku No. 1267 (2002), 112--141.


\bibitem[Me] {Me} {\sc J. P. Mestre},
  Construction de courbes de genre 2 \'a partir de leurs modules. In T.
  Mora and C. Traverso, editors,
 {\it Effective methods in algebraic geometry}, volume 94. {\it Prog.
 Math.}, 313-334. Birkh\"auser, 1991. Proc. Congress in Livorno, Italy,
April 17-21, (1990).

\bibitem[Sh] {Sh0} {\sc T. Shaska}, Curves of Genus Two Covering Elliptic
Curves, {\it PhD thesis}, University of Florida, 2001.


\bibitem[Sh1] {Sh1} {\sc T. Shaska},
Curves of genus 2 with $(n, n)$-decomposable Jacobians, J. Symbolic Comput. 31 (2001), no. 5, 603--617.


\bibitem[Sh2] {Sh2}
{\sc T. Shaska},
 Genus 2 curves with degree 3 elliptic subcovers,
 {\it Forum. Math.}, 2002, (accepted May 2002).


\bibitem[Sh3] {Sh3}
{\sc T. Shaska}, Genus 2 curves with (3,3)-split Jacobian and large automorphism group, Algorithmic Number Theory
(Sydney, 2002), {\bf 6}, 205-218, Lect. Not. in Comp. Sci., 2369, Springer, Berlin, 2002.

\bibitem[Sh4] {Sh4}
{\sc T. Shaska},
  Computational aspects of hyperelliptic curves,
 Computer Mathematics
(Beijing, 2003), Lect. Not. Ser. Comput., {\bf 10}, World Sci. Publishing, River Edge, NJ, 2003.


\bibitem[Sh5] {Sh5} {\sc T. Shaska}, Determining the automorphism group
of hyperelliptic curves, {\it Proceedings of the 2003 International Symposium on Symbolic and Algebraic
Computation}, 2003.


\bibitem[Sh6] {Sh6} {\sc T. Shaska}, Some special families of hyperelliptic curves,
 {\it J. Algebra Appl.}, 2003.

\bibitem[SV1] {ShV1}
{\sc T, Shaska and H. V\"olklein}, Elliptic subfields and automorphisms of genus two fields, {\it Algebra,
Arithmetic and Geometry with Applications.
 Papers from Shreeram S. Abhyankar's 70th Birthday Conference},
Springer (2003).

\bibitem[SV2] {ShV2}
{\sc T, Shaska and H. V\"olklein}, Genus two curves with degree 5 elliptic subcovers (preprint).

\bibitem[Shi] {Shi2} {\sc T. Shioda},
Constructing curves with high rank via symmetry. Amer. J. Math. 120 (1998), no. 3, 551--566.

\bibitem[Wi] {Wi} {\sc A. Wiman},
\"Uber die hyperelliptischen Curven vom den Geschlechte $p=4,5$, und 6, welche eindeutige Transformationen in sich
besitzen, {\it Bihang Kongl. Svenska Vetenskaps-Akademiens Handlingar} (1895), no. 21 (3), 1--41.

\end{thebibliography}
\end{document}